\documentclass[psamsfonts]{amsart}

\usepackage{amsfonts}

\usepackage[dvips]{graphicx}

\usepackage{amssymb}
\usepackage{amsmath}
\usepackage{color}

\usepackage[unicode,bookmarks,colorlinks]{hyperref}
\hypersetup{
    linkcolor=brickred,
}

\usepackage{amsmath,amstext,amssymb,amsopn,amsthm}
\usepackage{amsmath,amssymb,amsthm}
\usepackage[mathscr]{eucal}

\definecolor{mahogany}{cmyk}{0, 0.77, 0.87, 0}
\definecolor{salmon}{cmyk}{0, 0.53, 0.38, 0}
\definecolor{melon}{cmyk}{0, 0.46, 0.50, 0}
\definecolor{yellowgreen}{cmyk}{0.44, 0, 0.74, 0}
\definecolor{brickred}{cmyk}{0, 0.89, 0.94, 0.28}
\definecolor{OliveGreen}{cmyk}{0.64, 0, 0.95, 0.40}
\definecolor{RawSienna}{cmyk}{0, 0.72, 1.0, 0.45}
\definecolor{ZurichRed}{rgb}{1, 0, 0} 

\numberwithin{equation}{section}

\newtheorem{theorem}{Theorem}[section]

\newtheorem{prop}{Proposition}[section]
\newtheorem{corollary}[theorem]{Corollary}
\newtheorem{remark}[theorem]{Remark}

\newtheorem{definition}{Definition}

\newcommand{\R}{\mathbb{R}}

\newcommand{\Co}{\mathbb{C}}
\newcommand{\E}{\mathbb{E}}

\newcommand\F{\mathcal{F}}

\newcommand{\bb}[1]{\mathbb{#1}}
\newcommand{\bI}{\bb{I}}

\newcommand{\cadlag}{c\`adl\`ag}
\newcommand{\ds}{\displaystyle}
\newcommand{\g}{\mathfrak{g}}

\newcommand{\tr}{\mbox{tr}}

\newcommand{\Dom}{\mbox{Dom}}

\newcommand{\bean}{\begin{eqnarray*}}
\newcommand{\eean}{\end{eqnarray*}}

\newcommand{\calF}{\mathcal{F}}
\newcommand{\calB}{\mathcal{B}}
\newcommand{\calL}{\mathcal{L}}
\newcommand{\wG}{\widehat{G}}
\newcommand{\calS}{\mathcal{S}}

\newcommand\e{\varepsilon}

\newcommand{\bR}{\mathrm{I\! R\!}}

\date{}

\begin{document}

\title[]{Martingale transform and L\'evy Processes on Lie Groups}
\author{David Applebaum}\thanks{}
\address{School of Mathematics and Statistics, University of Sheffield,
Hicks Building, Hounsfield Road, Sheffield, England, S3 7RH}
\email{d.applebaum@sheffield.ac.uk}
\author{Rodrigo Ba\~nuelos}\thanks{R. Ba\~nuelos is supported in part  by NSF Grant
\# 0603701-DMS}
\address{Department of Mathematics, Purdue University, West Lafayette, IN 47907, USA}
\email{banuelos@math.purdue.edu}
\begin{abstract}  This paper constructs a class of martingale transforms based on L\'evy processes on Lie groups.  From these, a natural class of bounded linear operators on the $L^p$-spaces of the group (with respect to Haar measure) for $1<p<\infty$, are derived.  On compact groups these operators yield Fourier multipliers (in the Peter-Weyl sense) which include the second order Riesz transforms, imaginary powers of the Laplacian, and new classes of multipliers obtained by taking the L\'evy process to have conjugate invariant laws. Multipliers associated to subordination of the Brownian motion on the group are special cases of this last class.  These results extend (and the proofs simplify) those obtained in \cite{BanBieBog, BanBog} for the case of $\bR^n$.  An important feature of this work is the optimal nature of the $L^p$ bounds.

\end{abstract}

\maketitle
\tableofcontents

\section{Introduction}  Martingale inequalities have played  an important role in applications of probability to problems in analysis.  Many of these applications rest on the celebrated Burkholder-Davis-Gundy inequalities (commonly referred to these days as the {\it ``BDG inequalities"})  which compare the $L^p$-norms of martingales to the $L^p$-norms of their quadratic variations.   For a survey describing some of the early applications of these inequalities to several longstanding problems in analysis, we refer the reader to \cite{BanDav} and the many references therein. A predecessor to the  BDG inequalities which also has had a huge number of applications in analysis is Burkholder's 1966 inequality which proves the $L^p$ boundedness of martingale transforms \cite{Bur1}.  In 1984 Burkholder extended the martingale transform inequality in various directions  (including to the setting of Banach spaces with the unconditional martingale difference sequence property, the so called {\it ``UMD spaces"}) and obtained the sharp constants \cite{Bur2}. In recent years the sharp martingale transform inequalities have been widely applied to prove $L^p$-boundedness of Fourier multipliers on $\bR^n$, and in particular to the study of multipliers arising from the Riesz transforms, the Beurling-Ahlfors operator, operators of Laplace transform-type,  and related singular integrals. We refer the reader to \cite{Ban1} for references to this now very large literature and to some of the applications of these techniques to other problems of interest. These include connections to a longstanding conjecture of T. Iwaniec on the $L^p$-norm of Beurling-Ahlfors operator and a celebrated open problem of Morrey concerning rank-one convex and quasiconvex functions \cite[\S5]{Ban1}.  In \cite{BanBieBog} and  \cite{BanBog} sharp martingale transform inequalities were used to study a new class of Fourier multipliers which arise by certain natural transformations of L\'evy martingales associated to {\it L\'{e}vy processes}, i.e. the most general class of stochastic processes enjoying stationary and independent increments (of which Brownian motion is a special case.) A novel aspect of these new multipliers is that they are associated to processes whose paths contain  jump discontinuities and whose infinitesimal generators are non-local operators. A distinctive feature of Fourier multiplier inequalities obtained from sharp martingale transform inequalities is that in many instances the inequalities inherit the sharpness of the martingale transform inequalities. This happens, for example, in the important case of second order Riesz transforms and other multipliers with certain homogeneity properties \cite{BanOse}, including some of the  classical Marcinkiewicz-type multipliers \cite[p.110]{Ste2}.

The purpose of this paper is to show that the construction of the operators in \cite{BanBieBog} and  \cite{BanBog} can be carried out in the general setting of Lie groups, leading to a new collection of linear operators which are bounded on the $L^p$-spaces  of the group with bounds that are optimal with respect to  $p$ and which generalize in a natural way the classical Riesz transforms and the imaginary powers of the Laplacian. The results in \cite{BanBog} have recently seen applications to $L^p$ regularity for solutions to non-local elliptic problems in Euclidean space; see \cite{DonKim1, DonKim2}.  In the same way, we expect that the results of this paper will be of interest to those working in the analysis of Lie groups and its applications.

The paper is organized as follows.  In \S2.1, we recall the sharp Burkholder inequalities for martingales under the assumption of subordination.  While the discrete martingale versions serve as motivation, the results in this paper rest on versions of Burkholder's inequalities for continuous time martingales.   In \S2.2, we review some of the basic facts about L\'{e}vy processes on Lie groups and recall the fundamental theorem of Hunt which describes the structure of their infinitesimal generators.  In \S3, we define the martingale transforms and prove that they are bounded operators on $L^p$, for $1<p<\infty$ (Theorem \ref{mart-trasforms-bound}).  The martingale transforms lead naturally to operators defined on the Lie group with $L^p$ bounds which are exactly the same as those for martingale transforms (Corollary \ref{proj-transform-T} and Theorem \ref{projections1}).  In fact, the operators on the Lie group are limits of conditional expectations of the martingale transforms and for this reasons it is natural to call them ``projections of martingale transforms".  (For the latter point of view, we refer the reader to \cite{BanMen} and \cite[\S3.10]{Ban1} where this is done in the setting of $\bR^n$.)  In \S4 we specialize our construction to compact Lie groups where via the (non-commutative) Fourier transform we show that in various cases, including second order Riesz transforms and operators of Laplace transform-type,  our operators are Fourier multipliers. We show that a large class of multipliers can be obtained from central L\'{e}vy processes, i.e. those that have conjugate invariant laws. We obtain a new L\'{e}vy-Khintchine type formula for such processes on compact, connected semi-simple Lie groups which was anticipated in work of Liao \cite{Liao1, Liao}.  This formula plays an important role in the construction of multipliers but may also be of independent interest.  In \S5 we revisit, as a special case of the results for general Lie groups, the setting of $\bR^n$ and derive  the results of these two papers, rather directly, and without the technicalities of \cite{BanBog} and \cite{BanBieBog} concerning L\'evy systems.   But even more, our construction yields a larger class of operators than those studied in \cite{BanBog} and \cite{BanBieBog} (Corollary \ref{projection-Rn}). When further specializing the transformations, these operators consequently yield a wider class of Fourier multipliers \eqref{multiplier-3} than those in \cite{BanBog} and \cite{BanBieBog}.  Indeed the results of these papers are given by the multipliers in \eqref{multiplier-4} which are a sub-class of those presented in \eqref{multiplier-3}.

\vspace{5pt}

{\bf Notation.} If $G$ is a Lie group, ${\calB}(G)$ is the Borel $\sigma$-algebra of $G, C_{0}(G)$ is the Banach space of real-valued continuous functions on $G$ that vanish at infinity (equipped with the usual supremum norm), and $C_{c}^{\infty}(G)$ is the dense linear subspace of $C_{0}(G)$ comprising infinitely differentiable functions of compact support.  In any metric space $(M,d)$, the open ball of radius $r > 0$ centered on $p \in M$ is denoted by $B_{r}(p)$. If $a, b \in \R$, we define $a \vee b: = \max\{a,b\}$ and $a \wedge b:= \min\{a, b\}$. If $T > 0$ and $f:[0,T] \rightarrow \R$ has a left limit at $0 < s \leq T$ we write $f(s-):=\lim_{u \uparrow s}f(u)$ and we define the jump in $f$ at $s$ by $\Delta f(s) = f(s) - f(s-)$. Functions that are right continuous on $[0, T)$ and have left limits at every point in $(0,T]$ are said to be \cadlag~(from the French ``continue \`{a} droite et limit\'{e} \`{a} gauche''.)  For such functions the set $\{s \in (0, T); \Delta f(s) \neq 0\}$ is at most countable. Consequently we obtain identities such as $\int_{0}^{T}|f(s-)|^{p}ds = \int_{0}^{T}|f(s)|^{p}ds$ (for $1 \leq p < \infty$) which we will use freely within this paper. $M_n(\bR)$ denotes the space of all $n\times n$ matrices with real entries while $M_n(\Co)$ denotes the space of all $n \times n$ matrices with complex entries.  If $(M, \mathcal{M}, \mu)$ is a $\sigma$-finite measure space, the norm of $f \in L^{\infty}(M, \mathcal{M}, \mu; M_{n}(\R))$ is $||f||:= \mbox{ess.sup}_{x \in M}\sup\{|f(x)v|;v\in \bR^n, |v|\leq 1\}$. The norm in $L^{\infty}(M, \mathcal{M}, \mu; M_{n}(\Co))$ is defined similarly.

\section{Preliminaries}

\subsection{Sharp martingale inequalities}  We begin by recalling the celebrated sharp martingale transform inequalities of Burkholder.
 Let  $f=\{f_n, n\geq 0\}$ be a martingale (defined on some probability space) with difference sequence $d=\{d_k, k\geq 0\}$, where $d_k=f_k-f_{k-1}$ for $k\geq 1$ and $d_0=f_0$.  Given a predictable sequence of random variables $\{v_k, k\geq 0\}$ with $v_k\in [-1,1]$ a.s. for all $k$, the martingale difference sequence $\{v_kd_k, k\geq 0\}$ generates a new martingale called the {\it martingale transform} of $f$ and  denoted here by $g$.   We set $\|f\|_p=\sup_{n\geq 0}\|f_n\|_p$. In  \cite{Bur1} Burkholder proved that for all $1<p<\infty$, $\|g\|_p\leq C_p\|f\|_p$ for some constant $C_p$ depending only on $p$. In \cite{Bur2} he sharpened this result by proving  that under these assumptions  
 \begin{equation}\label{bur2}
\|g\|_p\leq (p^*-1)\|f\|_p,
\end{equation}
for all $1<p<\infty$, where $p^{*}:=\max\left\{p,\frac{p}{p-1}; 1 < p < \infty\right\}$ and that the constant $p^*-1$ is best possible.  In the sequel the constant $p^*-1$ will appear often.

Martingale inequalities of these type have a long history both in analysis and probability and we refer the reader to \cite{Ban1} for some of this literature and applications.
By considering dyadic martingales, inequality \eqref{bur2} contains the classical inequality of  Marcinkiewicz \cite{Mar} and Paley \cite{Pal} for Paley-Walsh martingales with the optimal constant.  That is, let $\{h_k,  k\geq 0\}$ be the Haar system in the Lebesgue unit interval $[0, 1)$ so that  $h_0=[0, 1), h_1=[0, \,1/2)-[1/2, \,1),  h_3=[0, \,1/4)-[1/4, \,1/2), h_4=[1/2, \,3/4)-(3/4,\, 1), \dots$, where the same notation is used for an interval as for its indicator function. Then for any sequence $\{a_k, {k\geq 0}\}$ of real numbers and any sequence $\{\e_k, {k\geq 0}\}$ of signs,
\begin{equation}\label{Paleyreal}
\Big\|\sum_{k=0}^{\infty} \varepsilon_k a_k h_k\Big\|_p\leq (p^*-1)\Big\|\sum_{k=0}^{\infty}  a_k h_k\Big\|_p, \qquad
\end{equation}
for $1<p<\infty$. The constant $p^*-1$ is best possible here as well.

An obvious question that arises from these sharp results is:  what happens to the best constant when the predictable sequence is non-symmetric in the sense that it takes values in $[0, 1]$ rather than in $[-1, 1]$?
In \cite{Cho}, K.P. Choi used the techniques of Burkholder to identify the best constant in this case as well.  While Choi's  constant is not as explicit (and simple) as the $p^*-1$ constant of Burkholder, one does have  considerable information about it. More precisely, Choi's result states that if  $v=\{v_k, k\geq 0\}$ with $v_k\in [0,1]$ a.s. for all $k$,  then
\begin{equation}\label{choi1}
\|g\|_p\leq c_p\|f\|_p,
\end{equation}
for all $1<p<\infty$ with the best constant $c_p$ satisfying
$$
c_p=\frac{p}{2}+ \frac{1}{2}\log\left(\frac{1+e^{-2}}{2}\right) +\frac{\alpha_2}{p}+\cdots
$$
where
$$\alpha_2=\left[\log\left(\frac{1+e^{-2}}{2}\right)\right]^2+\frac{1}{2}\log\left(\frac{1+e^{-2}}{2}\right)-2\left(\frac{e^{-2}}{1+e^{-2}}\right)^{2}. $$

Motivated by the inequalities of Burkholder  and Choi,  the following definition was introduced in \cite{BanOse}.

\begin{definition}\label{defC}
Let $-\infty<b<B<\infty$ and $1<p<\infty$ be given and fixed. We define $C_{p,b,B}$ as the least positive number $C$ such that for any real-valued martingale $f$ and for any transform $g$ of $f$ by a predictable sequence $v=\{v_k, k\geq 0\}$ with values almost surely in $[b,B]$,  we have
\begin{equation}\label{martin}
 ||g||_p\leq C||f||_p.
\end{equation}
\end{definition}

For any $0<a<\infty$, we see that $C_{p,-a,a}=a(p^*-1)$ by Burkholder's inequality \eqref{bur2} and
 $C_{p,0,a}=a\,c_p$ by Choi's inequality \eqref{choi1}.   An easy computation gives that for $b, B$ as above,

 \begin{equation}
 \max\left\{\left(\frac{B-b}{2}\right)(p^*-1), \, \max\{|B|, |b|\}\right\} \leq C_{p,b,B}\leq \max\{B, |b|\}(p^*-1).
  \end{equation}
However, the lack of any  general ``scaling" or ``translation" properties of $C_{p, b, B}$, outside the cases $[-a, a]$ and $[0, a]$, makes it very difficult to compute the constant.  For example, what is the value of $C_{p, 1, 2}$?

For the applications in this paper we require versions of the above inequalities for martingales indexed by continuous time. Suppose that $(\Omega,\mathcal{F},\mathbb{P})$ is a complete probability space, equipped with a right continuous filtration $(\mathcal{F}_t)_{t\geq 0}$  of sub-$\sigma$-algebras of $\mathcal{F}$.  We assume that $\F_0$ contains all the events of probability $0$. Let $X=\{X_t, t\geq 0\}$ and  $Y=\{Y_t, t\geq 0\}$ be adapted real valued martingales which have right-continuous paths with left-limits, i.e., {\cadlag}   martingales. We will denote by $[X,Y]$ the quadratic co-variation process of $X$ and $Y$.  When $X=Y$ we will simply write $[X]$ for $[X,X]$.  We refer the reader to Dellacherie and Meyer \cite{DelMey} for details on the construction and properties of $[X, Y]$ and $[X]$.
Following  \cite{BanWan} and  \cite{Wan}, we say that $Y$ is differentially subordinate to $X$  if $|Y_0| \leq |X_0|$ and $([X]_t-[Y]_t)_{t\geq 0}$ is nondecreasing and nonnegative as a function of $t$. We have the following extension of the Burkholder inequalities  proved in \cite{BanWan} for martingales with continuous paths (the so called ``Brownian martingales") and in \cite{Wan} in the general case.  More precisely,
if $Y$ is differentially subordinate to $X$, then
\begin{equation}\label{burkholder}
\|Y_T\|_p\leq (p^*-1)\|X_T\|_p,
\end{equation}
for all $1<p<\infty$ and all $T>0$.  The constant $p^*-1$ is best possible.

\begin{remark}\label{hilbert-valued}
We note here that Theorem \eqref{burkholder} holds in the exact same form if the martingales take values in a real or complex Hilbert space; see \cite{BanWan, Wan}. This is important for martingale transforms of complex valued functions such as the multipliers of Laplace transform-type given below.  That is, while we define martingale transforms with functions $A$ that take values in the set of all $n\times n$ matrices with real entries, $M_{n}(\R)$,  we could replace them with functions with values in $M_{n}(\Co)$.  A similar statement applies to the functions $\psi$ which we assume take values in $\bR$.  These changes to complex valued matrices $A$ and functions $\psi$ would not affect the bound $(p^*-1)$ given in the results that we obtain in this paper.
\end{remark}

We end this section by recalling a result for non-symmetric subordination proved in \cite{BanOse} which is the replacement for continuous time martingales of Choi's result.  This will be used below as well.
Suppose $-\infty<b<B<\infty$ and  $X_t$, $Y_t$ are two real valued martingales with right-continuous paths and  left-limits.  Suppose further that  $|Y_0|\leq |X_0|$ and that
\begin{equation}\label{banoseestiamte}
\left[\frac{B-b}{2}X,\frac{B-b}{2}X\right]_t-\left[Y-\frac{b+B}{2}X,Y-\frac{b+B}{2}X\right]_t\geq 0
\end{equation}
and is nondecreasing for all $t\geq 0$. (We call this property ``non-symmetric differential subordination".) Then
\begin{equation}\label{non-symm}
 ||Y_T||_p\leq C_{p,b,B}||X_T||_p, \quad 1<p<\infty,
\end{equation}
for all $T>0$.  The constant $C_{p,b,B}$ is best possible.

Note that when $B=a>0$ and $b=-a$, we have the case of \eqref{burkholder}.

\subsection{L\'evy processes on Lie groups}
Let $(\Omega, {\calF}, P)$
be a probability space and let $G$ be a Lie group of dimension $n$ with neutral element $e$ and Lie algebra $\g$. If $\phi: = (\phi(t), t \geq 0)$ is a stochastic process taking values in $G$ then the {\it right increment} of $\phi$ between $s$ and $t$ (where $s < t$) is the random variable $\phi(s)^{-1}\phi(t)$. We say that $\phi$ is a (left) {\it L\'{e}vy process} on $G$ if $\phi(0) = e$ (a.s.), $\phi$ has stationary and independent right increments and $\phi$ is stochastically continuous in that $\lim_{t \rightarrow 0} P(\phi(t) \in A) = 0$ for all $A \in \calB(G)$ for which $e \notin \overline{A}$.

Let $p_{t}(A):=P(\phi(t) \in A)$ for $A \in {\calB}(G), t \geq 0$ so that $p_{t}$ is the law of $\phi(t)$. Then $(p_{t}, t \geq 0)$ is a weakly continuous convolution semigroup of probability measures on $G$ for which $p_{0} = \delta_{e}$. Conversely given any such semigroup $(p_{t}, t \geq 0)$, we can always construct a L\'{e}vy process $(\phi_{t}, t \geq 0)$ on the space of all paths from $[0, \infty)$ to $G$ by using the celebrated Kolmogorov existence theorem.

For each $t \geq 0, f \in C_{0}(G), \sigma \in G$,  define $P_{t}f(\sigma) = \int_{G}f(\sigma \tau)p_{t}(d\tau)$. Then $(P_{t}, t \geq 0)$ is a (positivity preserving) $C_{0}$-contraction semigroup which commutes with left translations.  That is,  $L_{\sigma} P_{t} = P_{t} L_{\sigma}$ for all $t \geq 0, \sigma \in G$, where $L_{\sigma}f(\tau) = f(\sigma^{-1}\tau)$. In fact $(P_{t}, t \geq 0)$ is also a (almost everywhere positivity preserving) $C_{0}$-contraction semigroup on $L^{p}(G)$ for all $1 \leq p < \infty$ where $G$ is equipped with a {\it right-invariant} Haar measure (see e.g. \cite{Kun, App3}.)

We continue to work in $C_{0}(G)$ and let ${\calL}$ denote the infinitesimal generator of $(P_{t}, t \geq 0)$. So ${\calL}$ is a densely defined closed linear operator. We denote its domain by $\mbox{Dom}({\calL})$. We next give a precise description of ${\calL}$ that is due to Hunt \cite{Hunt} (see also \S3.1 in \cite{Liao}.)

Let $(X_{j}, 1 \leq j \leq n)$ be a fixed basis of $\g$ and define a dense linear manifold $C^{2}(G)$ in $C_{0}(G)$ by
$C^{2}(G): = \{ f \in C_{0}(G); X_{i}f
\in C_{0}(G)~ \mbox{and}~ X_{i}X_{j}f \in C_{0}(G) ~\mbox{for
all}~ 1 \leq i,j \leq n\}$. There exist functions ${x_{i} \in
C_{c}^{\infty}(G), 1 \leq i \leq n}$ so that $x_{i}(e) = 0,
X_{i}x_{j}(e) = \delta_{ij}$ and $(x_{1}, \ldots , x_{n})$ are canonical co-ordinates in a neighbourhood of $e$.

A measure $\nu$ defined on ${\calB}(G)$ is called a {\it
L\'{e}vy measure} whenever $\nu(\{e\}) = 0$,
\begin{equation}\label{levy-Lie}
\int_{U}\left(\sum_{i=1}^{n}x_{i}(\tau)^{2}\right)\nu(d\tau) < \infty\,\,\, \mbox{and}\,\,\, \nu(G-U) < \infty,
\end{equation}
for any Borel neighbourhood $U$ of $e$.

For the rest of the paper we will use the Einstein summation convention whereby we sum over repeated upper and lower indices.

\begin{theorem}[Hunt] \label{Hunt}  With the notation introduced above we have

\begin{enumerate}
\item $C^{2}(G) \subseteq \mbox{Dom}({\calL})$.
\bigskip
 \item For each
$\sigma \in G, f \in C^{2}(G)$,
\begin{eqnarray} \label{hu}
{\calL}f(\sigma) & = & b^{i}X_{i}f(\sigma) +
a^{ij}X_{i}X_{j}f(\sigma) \nonumber \\
  & + & \int_{G}(f(\sigma \tau) - f(\sigma) -
   x^{i}(\tau)X_{i}f(\sigma))\nu(d\tau),
\end{eqnarray}
where $b = (b^{1}, \ldots b^{n}) \in {\R}^{n}, a = (a^{ij})$ is a
non-negative-definite, symmetric $n \times n$ real-valued matrix
and $\nu$ is a L\'{e}vy measure on $G$.
\end{enumerate}
Conversely, any linear operator with a representation as above
is the restriction to $C^{2}(G)$ of the generator
corresponding to a unique convolution semigroup of probability
measures.
\end{theorem}

We call $(b, a, \nu)$ the characteristics of the L\'{e}vy process $\phi$.
From now on we will equip ${\calF}$ with a filtration $({\calF}_{t}, t \geq 0)$ of sub-$\sigma$-algebras and we will always assume that a given L\'{e}vy process $\phi$ is adapted to this filtration.  That is, $\phi(t)$ is ${\calF}_{t}$-measurable for each $t \geq 0$. We say that a L\'{e}vy process $\phi$ is \cadlag~if there exists $\overline{\Omega} \in {\calF}$ with $P(\overline{\Omega}) = 1$ such that the mappings $t \rightarrow \phi(t)(\omega)$ are right continuous with left limits existing for all $\omega \in \overline{\Omega}$.

Clearly any \cadlag~L\'{e}vy process $\phi$ is a Markov process with respect to its own filtration
and so for each $f \in C_{2}(G), t \geq 0$,

$$ f(\phi(t)) - f(e) - \int_{0}^{t}{\calL}f(\phi(s))ds$$
is a martingale
which we denote by $M_{f}(t)$ and which can be written as
\begin{equation}\label{decomposition}
M_{f}(t)= M_{f}^{c}(t) + M_{f}^{d}(t)
\end{equation}
where $M_{f}^{d}(t)$ and $M_{f}^{c}(t)$ are the discontinuous and continuous
parts, respectively, see \cite{DelMey}. In \cite{AK} these martingales were found to be
stochastic integrals against a Poisson random measure on $G$
and a Brownian motion in $\g$, respectively so that $\phi$ is in
fact the unique solution of the following stochastic differential equation

\begin{eqnarray} \label{apku}
f(\phi(t)) & = & f(e) + \int_{0}^{t}X_{i}f(\phi(s-))dB_{a}^{i}(s)
           + \int_{0}^{t}{\calL}f(\phi(s-))ds + \nonumber\\
&  & + \int_{0}^{t+}\int_{G}(f(\phi(s-)\sigma)
- f(\phi(s-))\tilde{N}(ds,d\sigma)
\end{eqnarray}
for each $f \in C_{2}(G), t \geq 0$. Here $B_{a} = (B_{a}(t), t \geq 0)$ is an
$n$-dimensional Brownian motion of mean zero and covariance matrix given by
Cov$(B_{a}^{i}(t)B_{a}^{j}(t)) = 2ta^{ij}$ for $t \geq 0, 1 \leq i,j \leq n$, and ${\tilde N}$
is the compensator defined for each $t\geq 0, E \in {\calB}(G)$ by

$$ {\tilde N}(t,E) = N(t,E) - t \nu(E),$$
where $N$ is a Poisson random measure on ${\bR}^{+} \times G$ with intensity
measure Leb$\times \nu$ which is independent of $B$. In fact,
$$ N(t, A) = \#\{0 \leq s \leq t, \Delta \phi(s) \in A\}$$
for each $t \geq 0, A \in {\calB}(G)$ where $\Delta \phi(s) = \phi(s-)^{-1}\phi(s)$ denotes the jump of the process at time $s > 0$.

Fix $T > 0$ and let $F: G \times G \rightarrow \R$ be continuous and $\xi \in L^{\infty}(\R^{+} \times G, \mbox{Leb}  \times \nu)$ be such that the mapping $s \rightarrow \xi(s, \tau)$ is left continuous for all $\tau \in G$  and assume that $\int_{0}^{T}\int_{G}\E(|F(\phi(s), \tau)|^{2})\nu(d \tau)ds < \infty$. We note that (see e.g. Lemma 4.2.2 in \cite{App1}, p.221) It\^{o}'s isometry in the framework of Poisson random measures yields,
\begin{eqnarray} \label{Itois}
\E\left(\left(\int_{0}^{T}F(\phi(s-), \tau)\xi(s, \tau)\tilde{N}(ds,d\tau)\right)^{2}\right) & = &
\int_{0}^{T}\int_{G}\E(|F(\phi(s), \tau)|^{2})|\xi(s, \tau)|^{2}\nu(d \tau)ds \nonumber \\
& \leq & ||\xi||^{2}\int_{0}^{T}\int_{G}\E(|F(\phi(s), \tau)|^{2})\nu(d \tau)ds,
\end{eqnarray}

In the sequel we will find it convenient to work with the standard Brownian motion $B = (B(t), t \geq 0)$ in $\R^{n}$ with covariance Cov$(B^{i}(t)B^{j}(t)) = t\delta^{ij}$. To implement this we choose an $n \times n$ matrix $\Lambda$ such that $\Lambda \Lambda^{T} = 2 a$ and define $Y_{i} \in \g$ by $Y_{i} = \Lambda_{i}^{j}X_{j}$ for $1 \leq i \leq n$. Then the integral with respect to Brownian motion in (\ref{apku}) may be rewritten
\bean
 \int_{0}^{t}X_{i}f(\phi(s-))dB_{a}^{i}(s) & = & \int_{0}^{t}Y_{i}f(\phi(s-))dB^{i}(s)\\
 & = & \int_{0}^{t}\nabla_{Y}f(\phi(s-)) \cdot dB(s),\eean
 where $\nabla_{Y}: = (Y_{1}, \ldots, Y_{n})$ and $\cdot$ is the usual inner product in $\R^{n}$.

 So far we have always assumed that the L\'{e}vy process $\phi$ starts at $e$ with probability one. In the sequel we will want to change the starting point to arbitrary $\rho \in G$ and we can achieve this by defining the \cadlag~Markov process $\phi^{(\rho)}(t) = \rho \phi(t)$ for each $t \geq 0$. Note that the process $\phi^{(\rho)}:=(\phi^{(\rho)}(t), t \geq 0)$ retains stationary and independent right increments. It follows easily from (\ref{apku}) that (with probability one)

\begin{eqnarray} \label{apku1}
f(\phi^{(\rho)}(t)) & = & f(\rho) +  \int_{0}^{t}\nabla_{Y}f(\phi^{(\rho)}(s-)) \cdot dB(s)
           + \int_{0}^{t}{\calL}f(\phi^{(\rho)}(s-))ds + \nonumber\\
&  & + \int_{0}^{t+}\int_{G}(f(\phi^{(\rho)}(s-)\sigma)
- f(\phi^{(\rho)}(s-))\tilde{N}(ds,d\sigma),
\end{eqnarray}

for each $f \in C_{2}(G), t \geq 0$.

\section{Martingale transforms and their projections}  In this section we will construct a family of operators which act on $L^p(G)$. These operators will arise as ``projections" of martingale transforms and will be bounded on $L^p(G)$ for all $1<p<\infty$.   To begin we fix $T > 0$ and choose $f \in C_{c}^{\infty}(G), \sigma \in G$. Then for all $0 \leq t \leq T$, we define

\begin{equation}\label{key}
M_{f}^{(T)}(\sigma, t): = (P_{T-t}f)(\sigma \phi(t))=(P_{T-t}f)(\phi^{\sigma}({t}))
\end{equation}
Later on we will also require the notation $M_{f}^{(T)}(\sigma):= M_{f}^{(T)}(\sigma, T)=f(\sigma\phi(T))$.

It follows from (\ref{apku1}) that
\begin{eqnarray}\label{mart1}
M_{f}^{(T)}(\sigma, t)&=&(P_{T}f)(\sigma)+ \int_{0}^{t}\nabla_{Y}(P_{T-s}f)(\phi^{(\sigma)}(s-)) \cdot dB(s)\nonumber\\
& +& \int_{0}^{t+}\int_{G}[(P_{T-s}f)(\phi^{(\sigma)}(s-)\tau)
- (P_{T-s}f)(\phi^{(\sigma)}(s-))]\tilde{N}(ds,d\tau),
\end{eqnarray}
and hence $M_{f,\sigma}^{(T)}:=(M_{f}^{(T)}(\sigma, t), 0 \leq t \leq T)$ is an $L^{2}$-martingale starting at $(P_{T}f)(\sigma)$. We fix a right-invariant Haar measure $m_{R}$ on $G$ and note that inside integrals we will always write $m_{R}(dg) = dg$.

\begin{definition}
Let $A \in L^{\infty}(\R^{+} \times G,  \mbox{Leb} \times m_{R};M_{n}(\R))$ and $\psi \in L^{\infty}(\R^{+} \times G \times G, \mbox{Leb} \times m_{R} \times \nu)$. We assume that $(A, \psi)$ are {\it regular} in that the mappings $(s,\rho) \rightarrow A(s, \rho)$ and $(s, \rho) \rightarrow \psi(s, \rho, \tau)$ (for all $\tau \in G$) are continuous and that
$||A|| \vee ||\psi|| \leq 1$. 
For each $f \in C_{c}^{\infty}(G), \sigma \in G, t \geq 0$ we set
\begin{eqnarray}\label{mart2}
M_{f}^{(T;A,\psi)}(\sigma, t) & = & \int_{0}^{t}A(T-s, \phi^{(\sigma)}(s-))\nabla_{Y}(P_{T-s}f)(\phi^{(\sigma)}(s-)) \cdot dB(s)\nonumber\\ & + & \int_{0}^{t+}\int_{G}\left\{(P_{T-s}f)(\phi^{(\sigma)}(s-)\tau)
- (P_{T-s}f)(\phi^{(\sigma)}(s-)\right\}\\
&\times&\left\{\psi(T-s, \phi^{(\sigma)}(s-), \tau)\right\} \tilde{N}(ds,d\tau)\nonumber.
\end{eqnarray}
This gives the new martingale
$$M_{f, \sigma}^{(T; A,\psi)}:=(M_{f}^{(T;A,\psi)}(\sigma, t), 0 \leq t \leq T)$$
which we shall call the martingale transform of  $M_{f, \sigma}^{(T)}$ by $(A, \psi)$.
\end{definition}

Computing the quadratic variations of both $(M_{f}^{(T)}(\sigma, t)$ and its transform $M_{f}^{(T;A,\psi)}(\sigma, t)$ (see e.g. equation (4.16) in \cite{App1} p.257) we find that
\begin{eqnarray}\label{qv1}
[(M_{f}^{(T)}(\sigma, \cdot)]_t & = & \int_{0}^{t}|\nabla_{Y}(P_{T-s}f)(\phi^{(\sigma)}(s-))|^{2}ds \nonumber \\ & + &
\int_{0}^{t+}\int_{G}[(P_{T-s}f)(\phi^{(\sigma)}(s-)\tau)
- (P_{T-s}f)(\phi^{(\sigma)}(s-))]^{2}N(ds,d\tau)
\end{eqnarray}
while
\begin{eqnarray}\label{qv2}
[(M_{f}^{(T;A, \psi)}(\sigma, \cdot)]_t & = & \int_{0}^{t}|A(T-s, \phi^{(\sigma)}(s-))\nabla_{Y}(P_{T-s}f)(\phi^{(\sigma)}(s-))|^{2}ds \nonumber \\ & + &
\int_{0}^{t+}\int_{G}\left\{[(P_{T-s}f)(\phi^{(\sigma)}(s-)\tau)
- (P_{T-s}f)(\phi^{(\sigma)}(s-))]^{2}\right\}\nonumber\\&\times&\left\{[\psi(T-s, \phi^{(\sigma)}(s-), \tau)]^{2}\right\}N(ds,d\tau).
\end{eqnarray}

From these formulas, the fact that $0=|M_{f}^{(T;A,\psi)}(\sigma, 0)|\leq |M_{f}^{(T)}(\sigma, 0)|$, and our assumption that $||A|| \vee ||\psi|| \leq 1$, it follows  that $M_{f, \sigma}^{(T; A,\psi)}$ is differentially subordinate to $M_{f,\sigma}^{(T)}$.  Now, for $1<p<\infty$, set
$$||X||_{p} = \left(\int_{G}\E(|X(\sigma)|^{p})d\sigma)\right)^{\frac{1}{p}},$$
 for $X \in L^{p}(\Omega \times G)$.  Applying Burkholder's result in the form given by the inequality  \eqref{burkholder}, we obtain

\begin{equation}\label{fundineq}
||M_{f}^{(T;A,\psi)}||_{p} \leq (p^{*}-1)||M_{f}^{(T)}||_{p}.
\end{equation}

But using Fubini's theorem and the right invariance of the Haar measure we see that
\begin{eqnarray}\label{ineq1}
||M_{f}^{(T)}||_{p}^p&=&||f(\cdot \phi(T)||_{p}^{p}\nonumber\\
 & = & \int_{G}\int_{G}|f(\sigma \tau)|^{p}p_{T}(d\tau)d\sigma\\
& = & \int_{G}|f(\sigma)|^{p}d\sigma = ||f||_{p}^{p}\nonumber.
\end{eqnarray}
We summarize the above calculations in the following

\begin{theorem}\label{mart-trasforms-bound} Let $A \in L^{\infty}(\R^{+} \times G,  \mbox{Leb} \times m_{R};M_{n}(\R))$ and $\psi \in L^{\infty}(\R^{+} \times G \times G, \mbox{Leb} \times m_{R} \times \nu)$ with $(A,\psi)$ regular and
$||A|| \vee ||\psi|| \leq 1$.  For any $0<T<\infty$, the map $f\to M_{f}^{(T;A,\psi)}$ defines a linear operator from $L^p(G)\to L^{p}(\Omega \times G)$, $1<p<\infty$,  with
\begin{equation}\label{fundineq2}
\|M_{f}^{(T;A,\psi)}\|_{p} \leq (p^{*}-1)\|f\|_{p},
\end{equation}
for $1<p<\infty$. In particular, the bound is independent of $T$.

\end{theorem}

Now let $q = \frac{p}{p-1}$ and for given $g \in C_{c}^{\infty}(G)$, we define a linear functional $\Lambda_{g}^{T;A,\psi}$ on $C_{c}^{\infty}(G)$ by the prescription
\begin{equation}\label{lf}
\Lambda_{g}^{T;A,\psi}(f) = \int_{G}\E(M_{f}^{(T;A,\psi)}(\sigma)M_{g}^{(T)}(\sigma))d\sigma
\end{equation}

Using H\"{o}lder's inequality, inequality (\ref{fundineq}) and equality (\ref{ineq1}), we obtain
\begin{eqnarray}\label{ineq2}
|\Lambda_{g}^{T;A,\psi}(f)| & \leq & ||M_{f}^{(T;A,\psi)}||_{p}||M_{g}^{(T)}||_{q}\nonumber\\
& \leq & (p^{*} -1)||M_{f}^{(T)}||_{p}||M_{g}^{(T)}||_{q}\nonumber\\
& \leq & (p^{*} -1)||f||_{p}||g||_{q}.
\end{eqnarray}

Hence $\Lambda_{g}^{T;A,\psi}$ extends to a bounded linear functional on $L^{p}(G)$ and by duality, there exists a bounded linear operator $S_{A,\psi}^{T}$ on $L^{p}(G)$  for which
$$ \Lambda_{g}^{T;A,\psi}(f) = \int_{G}S_{A,\psi}^{T}f(\sigma)g(\sigma)d\sigma, $$
for all $f \in L^{p}(G), g \in L^{q}(G)$ and with $||S_{A,\psi}^{T}||_{p} \leq (p^{*}-1)$.
We summarize this in the following

\begin{corollary}\label{proj-transform-T} Let $A \in L^{\infty}(\R^{+} \times G, \mbox{Leb} \times m_{R};M_{n}(\R))$ and $\psi \in L^{\infty}(\R^{+} \times G \times G, \mbox{Leb} \times m_{R} \times \nu)$ with $(A,\psi)$ regular and
$||A|| \vee ||\psi|| \leq 1$.  For any $0<T<\infty$, the map $f\to S_{A,\psi}^{T}f$ defines a linear operator from $L^p(G)\to L^{p}(G)$, $1<p<\infty$,  with
\begin{equation}\label{fundineq3}
\|S_{A,\psi}^{T}f\|_{p} \leq (p^{*}-1)\|f\|_{p},
\end{equation}
for $1<p<\infty$.  In particular, the bound is independent of $T$.
\end{corollary}
Next we probe the structure of (\ref{lf}) using It\^{o}'s isometry for stochastic integrals driven by both Brownian motion and a Poisson random measure (see \eqref{Itois} for the latter) with the aim of letting $T\to\infty$ in $S_{A,\psi}^{T}$. Again, using Fubini's theorem and right invariance of Haar measure as in the derivation of (\ref{ineq1}), we obtain

\begin{eqnarray}\label{PL}
\Lambda_{g}^{T;A,\psi}(f)&=&  \int_{G}\int_{0}^{T}\E\{A(T-s, \phi^{(\sigma)}(s-))\nabla_{Y}(P_{T-s}f)(\phi^{(\sigma)}(s-))\nonumber\\
&\times&\cdot \nabla_{Y}(P_{T-s}g)(\phi^{(\sigma)}(s-))\}dsd\sigma \nonumber \\
& + & \int_{G}\int_{G}\int_{0}^{T}\E\{[(P_{T-s}f)(\phi^{(\sigma)}(s-)\tau)
- (P_{T-s}f)(\phi^{(\sigma)}(s-))] \nonumber \\ & \times & [(P_{T-s}g)(\phi^{(\sigma)}(s-)\tau)
- (P_{T-s}g)(\phi^{(\sigma)}(s-))]\nonumber \\
&\times&\psi(T-s, \phi^{(\sigma)}(s-),\tau)\}ds\nu(d\tau)d\sigma \nonumber \\
& = & \int_{0}^{T}\int_{G} A(T-s, \sigma)\nabla_{Y}(P_{T-s}f)(\sigma)\cdot \nabla_{Y}(P_{T-s}g)(\sigma)d\sigma ds \nonumber \\
& +&  \int_{0}^{T}\int_{G}\int_{G}[(P_{T-s}f)(\sigma\tau)
- (P_{T-s}f)(\sigma)][(P_{T-s}g)(\sigma \tau)
- (P_{T-s}g)(\sigma)]\nonumber\\
&\times& \psi(T-s, \sigma,\tau)\nu(d\tau)d\sigma ds \nonumber \\
& = & \int_{0}^{T}\int_{G} A(s, \sigma)\nabla_{Y}(P_{s}f)(\sigma)\cdot \nabla_{Y}(P_{s}g)(\sigma)d\sigma ds\nonumber\\
& + & \int_{0}^{T}\int_{G}\int_{G}[(P_{s}f)(\sigma\tau)
- (P_{s}f)(\sigma)][(P_{s}g)(\sigma \tau)
- (P_{s}g)(\sigma)]\nonumber \\
&\times&
\psi(s,\sigma,\tau)\nu(d\tau)d\sigma ds
\end{eqnarray}

Choosing
$$ A(s, \sigma) = \ds\frac{\nabla_{Y}(P_{s}f)(\sigma) \otimes \nabla_{Y}(P_{s}g)(\sigma)}{|\nabla_{Y}(P_{s}f)(\sigma)|\,|\nabla_{Y}(P_{s}g)(\sigma)|},$$
for $s \geq 0, \sigma \in G$ and observing that $\|A\|=1$, inequality \eqref{ineq2} (with $\psi\equiv 0$) gives that
\begin{equation}\label{LP1}
\int_{0}^{T}\int_{G} |\nabla_{Y}(P_{s}f)(\sigma)|\,|\nabla_{Y}(P_{s}g)(\sigma)|d\sigma ds\leq (p^{*} -1)||f||_{p}||g||_{q}.
\end{equation}
Thus  the first integral in \eqref{PL} converges absolutely as $T \rightarrow \infty$ and
\begin{equation}\label{LP2}
\int_{0}^{\infty}\int_{G} |\nabla_{Y}(P_{s}f)(\sigma)|\,|\nabla_{Y}(P_{s}g)(\sigma)|d\sigma ds\leq (p^{*} -1)||f||_{p}||g||_{q}.
\end{equation}

We now consider the second integral  in \eqref{PL}.  This time we take $A\equiv 0$ and
$$
\psi(s, \sigma, \tau) = \mbox{sign}((P_{s}f)(\sigma \tau) - P_{s}f(\sigma))((P_{s}g)(\sigma \tau) - P_{s}g(\sigma)))
$$
and obtain that
\begin{equation*}\label{secondterm}
\int_{0}^{T}\int_{G}\int_{G}|(P_{s}f)(\sigma \tau) - P_{s}f(\sigma))| |(P_{s}g)(\sigma \tau) - P_{s}g(\sigma))|\nu(d\tau)d\sigma ds
\leq (p^{*} -1)||f||_{p}||g||_{q},
\end{equation*}
with the right hand side independent of $T$ and hence we can also let $T\to\infty$ in the second integral of \eqref{PL}.

We summarize the above result  in the following

\begin{theorem}\label{projections1} Let $A \in L^{\infty}(\R^{+} \times G,\mbox{Leb} \times m_{R};M_{n}(\R))$ and $\psi \in L^{\infty}(\R^{+} \times G \times G, \mbox{Leb} \times m_{R} \times \nu)$ with  $(A, \psi)$ regular and
$||A|| \vee ||\psi|| \leq 1$. There exists a bounded linear operator $S_{A, \psi}$ on $L^{p}(G)$, $1<p<\infty$, for which

\begin{eqnarray} \label{identity}
\int_{G}S_{A,\psi}f(\sigma)g(\sigma)d\sigma & = & \int_{0}^{\infty}\int_{G} A(s, \sigma)\nabla_{Y}(P_{s}f)(\sigma)\cdot \nabla_{Y}(P_{s}g)(\sigma)d\sigma ds \nonumber \\
& + & \int_{0}^{\infty}\int_{G}\int_{G}[(P_{s}f)(\sigma\tau)
- (P_{s}f)(\sigma)][(P_{s}g)(\sigma \tau)
- (P_{s}g)(\sigma)]\\
&\times& \psi(s,\sigma,\tau)\nu(d\tau)d\sigma ds,\nonumber
\end{eqnarray}
for al $f, g \in C_{c}^{\infty}(G)$.  Furthermore, for all $f\in L^p(G)$ and $g\in L^q(G)$, $\frac{1}{p}+\frac{1}{q}=1$,
\begin{equation}\label{bound1}
\Big|\int_{G}S_{A,\psi}f(\sigma)g(\sigma)d\sigma \Big| \leq (p^*-1)\|f\|_p\,\|g\|_q
\end{equation}
and
\begin{equation}\label{bound2}
\|S_{A,\psi}f\|_p\leq (p^*-1)\|f\|_p.
\end{equation}

\end{theorem}

Suppose $A \in L^{\infty}(\R^{+} \times G,\mbox{Leb} \times m_{R};M_{n}(\R))$ is symmetric with the property that for all $\xi\in \R^n$,
\begin{equation}\label{non-symm-A}
b|\xi|^2\leq A(s, \sigma)\xi\cdot \xi \leq B|\xi|^2
\end{equation}
for all $(s, \sigma)\in \R^{+} \times G$, $-\infty<b<B<\infty$.  Then by a simple computation (see \cite{BanOse}, \S4)  the martingale $Y=M_{f, \sigma}^{(T; A, 0)}$
is subordinate to the martingale  $X=M_{f, \sigma}^{(T)}$  in the sense of inequality  \eqref{banoseestiamte}.  From the construction of our operators $S_{A, \psi}$ we have (in this case the assumption $||A|| \leq 1$ is no longer needed):
\begin{theorem}\label{non-symm-thm} Let
$A \in L^{\infty}(\R^{+} \times G, \mbox{Leb} \times m_{R};M_{n}(\R))$ be continuous, symmetric and satisfying \eqref{non-symm-A}.  Then for all $1<p<\infty$ and $f\in L^p(G)$,
\begin{equation}\label{nonsymmetric1}
\|S_{A, 0}f\|_p\leq C_{p,b,B}\|f\|_p,
\end{equation}
where $C_{p, b, B}$ is the constant in the inequality \eqref{non-symm}
\end{theorem}

As we shall see below (see \S5), the inequalities \eqref{bound2} and \eqref{nonsymmetric1} are sharp as they include the bounds for Riesz transforms in $\bR^n$ and other homogeneous multipliers (see \cite{BanOse}).

\section{Compact groups:  Fourier multipliers}

In this action we are interested in cases of a compact Lie group where the operators $S_{A,\psi}$ are Fourier multipliers. Thus from now on we will assume that $G$ is a compact Lie group. Then every (left or right) Haar measure on $G$ is bi-invariant and also finite. In the sequel we will always assume that Haar measure is normalized to have total mass one. Let $\wG$ be the {\it unitary dual} of $G$, i.e. the set of all equivalence classes (modulo unitary equivalence) of irreducible representations of $G$. In the sequel we will often identify a class in $\wG$ with a representative element. The set $\wG$ is countable and each $\pi \in \wG$ is finite dimensional, so that that there exists a finite dimensional complex Hilbert space $V_{\pi}$, having dimension $d_{\pi}$ such that for each $\sigma \in G, \pi(\sigma)$ is a unitary matrix. We define the co-ordinate functions $\pi_{ij}(\sigma) = \pi(\sigma)_{ij}$ for $\sigma \in G, 1 \leq i,j \leq d_{\pi}$. We have $\pi_{ij} \in C^{\infty}(G)$ and the celebrated Peter-Weyl theorem tells us that $\{\sqrt{d_{\pi}}\pi_{ij}; 1 \leq i,j \leq d_{\pi}, \pi \in \wG\}$ is a complete orthonormal basis for $L^{2}(G, \Co)$. For each $f \in L^{2}(G, \Co)$, we define its {\it non-commutative Fourier transform} to be the matrix $\widehat{f}(\pi)$ defined by
$$ \widehat{f}(\pi)_{ij} = \int_{G}f(\sigma)\pi_{ij}(\sigma^{-1})d\sigma,$$ for each $1 \leq i,j \leq d_{\pi}$.
We will have need of the following version of the {\it Plancherel theorem} (which can be found in e.g. \cite{Far}, Theorem 6.4.2, p.110)
\begin{equation} \label{Planch}
\int_{G}f(\sigma)\overline{g(\sigma)}d\sigma = \sum_{\pi \in \wG}d_{\pi}\tr(\widehat{f}(\pi)\widehat{g}(\pi)^{*}),
\end{equation}
for $f, g\in L^2(G,  \Co)$.  In particular if $T$ is a bounded linear operator on $L^{2}(G, \Co)$ we have
\begin{equation} \label{Planch1}
\int_{G}Tf(\sigma)\overline{g(\sigma)}d\sigma = \sum_{\pi \in \wG}d_{\pi}\tr(\widehat{Tf}(\pi)\widehat{g}(\pi)^{*}).
\end{equation}

We say that the operator $T$ is a {\it Fourier multiplier} if for each $\pi \in \wG$ there exists a $d_{\pi} \times d_{\pi}$ complex matrix $m_{T}(\pi)$ so that

\begin{equation} \label{FM}
\widehat{Tf}(\pi) = m_{T}(\pi)\widehat{f}(\pi).
\end{equation}
We call the matrices $(m_{T}(\pi), \pi \in \wG)$ the {\it symbol} of the operator $T$.\footnote{See \cite{RT} for a monograph treatment of pseudo differential operators and their symbols on compact Lie groups and \cite{RW} for a study of Fourier multipliers from this perspective. For probabilistic developments in the spirit of the present work, see \cite{App, App2}.}

Given $\pi \in \wG$ we obtain the {\it derived representation} $d\pi$ of the Lie algebra $\g$ from the identity
$$ \pi(\exp(X)) = e^{d\pi(X)}, $$
for each $X \in \g$. Then $d\pi(X)$ is a skew-hermitian matrix acting in $V_{\pi}$. We now equip $\g$ with an Ad-invariant metric (which induces a bi-invariant Riemannian metric on $G$.) From now on $\{X_{1}, \dots, X_{n}\}$ will be an orthonormal basis for $\g$ with respect to the given metric. We define the {\it Casimir operator} $\Omega_{\pi}$ by $\Omega_{\pi}: = \sum_{i=1}^{n}d\pi(X_{i})^{2}$. Then it can be shown (see e.g. \cite{Far} Corollary 6.7.2, p.122) that $\Omega_{\pi} = -\kappa_{\pi}I_{\pi}$ where $I_{\pi}$ is the identity matrix acting on $V_{\pi}$ and $\kappa_{\pi} \geq 0$ (with $\kappa_{\pi}=0$ if and only if $\pi$ is the trivial representation).
The Laplace-Beltrami operator $\Delta = \sum_{i=1}^{n}X_{i}^{2}$ is an essentially self-adjoint operator in $L^{2}(G, \Co)$ with domain $C(G, \Co)$ having discrete spectrum with
\begin{equation} \label{spec}
\Delta \pi_{ij} = - \kappa_{\pi}\pi_{ij},
\end{equation}
for all $\pi \in \wG, 1 \leq i,j \leq d_{\pi}$.

\subsection{Brownian motion}

Brownian motion (with twice the usual auto-covariance) on a Lie group $G$ is the L\'{e}vy process having characteristics $(b, I, 0)$. It is well known (see e.g. \cite{Elw}, \S6) that it can be obtained as the unique solution of the Stratonovich stochastic differential equation
\begin{equation} \label{sde}
d\phi(t) = \sqrt{2}X^{i}(\phi(t)) \circ dB_{i}(t),
\end{equation}  with initial condition $\phi(0) = e$ (a.s.) In this case $(P_{t}, t \geq 0)$ is the heat semigroup with generator ${\calL} = \Delta$.

The following is well-known but we include a short proof for the reader's convenience:

\begin{prop} \label{ftheat} For all $t \geq 0, \pi \in \wG, f \in L^{2}(G, \Co)$
$$ \widehat{P_{t}f}(\pi) = e^{-t\kappa_{\pi}}I_{\pi}.$$
\end{prop}

{\it Proof.} It follows from (\ref{spec}) that for all $1 \leq i,j \leq d_{\pi}$,
     $$P_{t}\pi_{ij} = e^{-t\kappa_{\pi}}\pi_{ij}.$$
But then
\bean \widehat{P_{t}f}(\pi)_{ij} & = & \int_{G}P_{t}f(g)\pi(g^{-1})_{ij}dg\\
& = &  \int_{G}P_{t}f(g)\overline{\pi(g)_{ji}}dg\\
& = & \int_{G}f(g)\overline{P_{t}\pi(g)_{ji}}dg\\
& = & e^{-t\kappa_{\pi}}\int_{G}f(g)\pi(g^{-1})_{ij}dg \\
& = & e^{-t\kappa_{\pi}}\widehat{f}(\pi)_{ij}, \eean
as required. $\hfill \Box$

\subsection{Second order Riesz transforms}

The {\it first order Riesz transform} in the direction $X \in \g$ is the operator $R_{X} = X(-\Delta)^{-\frac{1}{2}}$. It is shown in \cite{Arc} (using the martingale inequalities from \cite{BanWan}) that if $G$ is endowed with a bi-invariant Riemannian metric and $|X| = 1$ then $R_{X}$ is a bounded operator on $L^{p}(G)$ and $||R_{X}||_{p} \leq B_{p}$ for all $1 < p < \infty$, where $B_{p}: = \cot\left(\frac{\pi}{2p^{*}}\right)$ is the ``Pichorides constant'' (i.e. the $L^{p}$ norm of the Hilbert transform on the real line.) We write $R_{j}: = R_{X_{j}}$, for $1 \leq j \leq n$, so that $R_{i}R_{j}=X_{i}X_{j}{\Delta}^{-1}$. Let $C$ be a $n \times n$ matrix with $||C|| \leq 1$. We define the {\it second order Riesz transform} (see also \cite{GMS, BanBau}) $R_{C}^{(2)}$ to be

$$ R_{C}^{(2)}: = \sum_{i, j = 1}^{n}C_{ji}R_{i}R_{j} = \sum_{i, j = 1}^{n}C_{ji}X_{i}X_{j}\,{\Delta}^{-1}. $$

We will show that $R_{C}^{(2)}$ is precisely an operator of the form $S_{A, \psi}$ where $\psi = 0$. First we show that $R_{C}^{(2)}$ is a Fourier multiplier. Indeed this follows by using (\ref{Planch1}) and computing for $C^{\infty}(G)$
\bean \int_{G}R_{C}^{(2)}f(\sigma) g(\sigma)d\sigma & = & \sum_{i, j = 1}C_{ji}\sum_{\pi \in \wG}d_{\pi}\tr(\widehat{R_{i}R_{j}f}(\pi)\widehat{g}(\pi)^{*})\\
& = & -  \sum_{i, j = 1}C_{ji}\sum_{\pi \in \wG}d_{\pi}\frac{1}{\kappa_{\pi}}\tr(d\pi(X_{i})d\pi(X_{j})\widehat{f}(\pi)\widehat{g}(\pi)^{*})\\
& = & \sum_{\pi \in \wG}d_{\pi}\tr(\Omega_{C}(\pi)\widehat{f}(\pi)\widehat{g}(\pi)^{*}), \eean
where for all $\pi \in \wG$,
$$ \Omega_{C}(\pi): = -\frac{1}{\kappa_{\pi}}\sum_{i, j = 1}^{n}C_{ji}d\pi(X_{i})d\pi(X_{j}). $$
Hence we conclude that $R_{C}^{(2)}$ is a Fourier multiplier with symbol $\Omega_{C}(\cdot)$.

We connect to the work of the previous section by taking $\psi = 0$ and $\phi$ to be a Brownian motion. We write $S_{A} = S_{A,0}$ and we take $A$ to be a constant matrix.
Using (\ref{identity}), (\ref{Planch1}) and Proposition \ref{ftheat}, we obtain, for $f, g \in C^{\infty}(G)$,
\bean \int_{G}S_{A}f(\sigma)g(\sigma)d\sigma & = & 2\int_{0}^{\infty}\int_{G} A \nabla_{X}(P_{s}f)(\sigma)\cdot \nabla_{X}(P_{s}g)(\sigma)d\sigma ds\\
& = &  2\sum_{i,j=1}^{n}A_{ij}\int_{0}^{\infty}\int_{G} X_{i}(P_{s}f)(\sigma)X_{j}(P_{s}g)(\sigma)d\sigma ds\\
& = &  2\sum_{i,j=1}^{n}A_{ij}\int_{0}^{\infty}\sum_{\pi \in \wG}d_{\pi}e^{-2s\kappa_{\pi}}\tr(d\pi(X_{i})\widehat{f}(\pi)\widehat{g}(\pi)^{*}(-d\pi(X_{j}))ds\\
& = & -\sum_{i,j=1}^{n}A_{ij}\sum_{\pi \in \wG}d_{\pi}\frac{1}{\kappa_{\pi}}\tr(d\pi(X_{j})d\pi(X_{i})\widehat{f}(\pi)\widehat{g}(\pi)^{*})\\
& = & \sum_{\pi \in \wG}d_{\pi}\tr(\Omega_{A}(\pi)\widehat{f}(\pi)\widehat{g}(\pi)^{*}), \eean
and so we deduce that $R_{C}^{(2)} = S_{C}$.

The following is a corollary of the above discussion and Theorem \ref{projections1}, inequality \eqref{bound1}.

\begin{corollary}  Consider the second order Riesz transform defined by
\begin{equation}\label{secondriesz1}
R_{C}^{(2)}: = \sum_{i, j = 1}^{n}C_{ji}R_{i}R_{j} = \sum_{i, j = 1}^{n}C_{ji}X_{i}X_{j}{\Delta}^{-1}
\end{equation}
with symbol
\begin{equation}\label{secondriesz2}
 \Omega_{C}(\pi): = -\frac{1}{\kappa_{\pi}}\sum_{i, j = 1}^{n}C_{ji}d\pi(X_{i})d\pi(X_{j}).
 \end{equation}
 Then for $1<p<\infty$, $f\in L^p(G)$,
 \begin{equation}
 \|R_{C}^{(2)}f\|_p\leq (p^*-1)\|f\|_p.
 \end{equation}
\end{corollary}

The next corollary follows from Theorem \ref{non-symm-thm}.  We record  it here in the form stated in \cite{BanOse} in the case of $\R^n$ (see also \cite{BanBau} for compact Lie groups). Note that we here drop the condition $||C|| \leq 1$ as in Theorem \ref{non-symm-thm}.
\begin{corollary}\label{Riesz1}
Let $C$ be an $n\times n$ symmetric matrix  with eigenvalues $\lambda_1\leq \lambda_2\leq \ldots \leq \lambda_n$.
Then for $1<p<\infty$, $f\in L^p(G)$,
\begin{equation}\label{generalriesz1}
 \|R_{C}^{(2)}f\|_p\leq C_{p,\,\lambda_1,\lambda_d}\|f\|_p,
\end{equation}
where $C_{p,\,\lambda_1,\lambda_d}$ is the constant in \eqref{martin}.   In particular, if $J\subsetneq \{1,2,\,\ldots,n\}$, then
\begin{equation}\label{squareriesz2}
\left|\left|\sum_{j\in J} R_j^2 f\right|\right|_p\leq C_{p,0,1}\|f\|_p=c_p\|f\|_p, \quad 1<p<\infty,
\end{equation}
where $c_p$ is the Choi constant in \eqref{choi1}.
\end{corollary}

We remark that in the case of $\bR^n$  both bounds above are sharp,  see \cite{BanOse} and  \cite{GMS}.

 \subsection{Multipliers of Laplace transform-type}

We continue to work with Brownian motion (as in the previous subsection) and this time we take $A(s):=A(s,\cdot)$ to be a time-dependent matrix-valued function. Then (\ref{identity}), integration by parts and (\ref{Planch1}) yields for all $f, g \in C^{\infty}(G)$,
\bean \int_{G}S_{A}f(\sigma)g(\sigma)d\sigma & = &2 \int_{0}^{\infty}\int_{G} A(s) \nabla_{X}(P_{s}f)(\sigma)\cdot \nabla_{X}(P_{s}g)(\sigma)d\sigma ds,\\
& = &-2 \int_{0}^{\infty}\int_{G}A(s)\Delta P_{2s}f(\sigma)g(\sigma)d\sigma ds.\\
& = & 2\int_{0}^{\infty}A(s)\sum_{\pi \in \wG}d_{\pi}\kappa_{\pi}e^{-2s \kappa_{\pi}}\tr(\widehat{f}(\pi)\widehat{g}(\pi))ds. \eean
Hence for all $\pi \in \wG$ we have

\begin{equation} \label{LT}
\widehat{S_{A}f}(\pi) = \left(\int_{0}^{\infty}2\kappa_{\pi}e^{-2s \kappa_{\pi}}A(s)ds\right)\widehat{f}(\pi)
\end{equation}
and $S_{A}$ is an operator of {\it Laplace transform-type} (see \cite{Stein} p.58, \cite{Ban1}, \S3.11). It is clearly a Fourier multiplier in the sense of (\ref{FM}). A special case of particular interest is obtained by taking $A(s) = \ds\frac{(2s)^{-i\gamma}}{\Gamma(1-i\gamma)}{I}$,  where $\gamma \in \R$. In this case an easy computation shows that the symbol is $m_{S_{A}}(\pi) = \kappa_{\pi}^{-i\gamma}$ and so $S_{A} = (-\Delta)^{i\gamma}$ is an imaginary power of the Laplacian (see also \cite{SW}) and we have the following corollary of Theorem \ref{projections1}, inequality \eqref{bound1}. (Note that for this case we need to apply the Burkholder inequality for complex valued martingales using the results in \cite{BanWan} or \cite{Wan}, as indicated in Remark \ref{hilbert-valued}).

\begin{corollary}  For $1<p<\infty$, $f\in L^p(G)$, we have
\begin{equation}
\|(-\Delta)^{i\gamma}f\|_{p} \leq \frac{p^{*}-1}{|\Gamma(1- i\gamma)|}\|f\|_{p}.
\end{equation}
\end{corollary}

This bound was first proved in $\bR^n$ in \cite[Eq. (7.3)]{Hyt}.  For further discussion on this result in $\bR^n$ and comparison to previous known bounds, see \cite{Ban1}.

\subsection{Subordinate Brownian motion}

In this section we will consider an operator of the form $S_{\psi}:=S_{0, \psi}$ which is built from the non-local part of the L\'{e}vy process. To this end, let $(T(t), t \geq 0)$ be a subordinator, i.e. an almost surely non-decreasing real-valued L\'{e}vy process (so $T(t)$ takes values in $\R^{+}$ with probability one for all $t \geq 0$) . Then we have (see e.g. \cite{SSV}) that $\E(e^{-uT(t)}) = e^{-t h(u)}$ for all $u > 0, t \geq 0$ where $h:(0, \infty) \rightarrow \R^{+}$ is a Bernstein function for which $\lim_{u \rightarrow 0}h(u) = 0$. Hence $h$ must take the form
\begin{equation} \label{bern}
      h(u) = c u + \int_{(0, \infty)}(1 - e^{-uy})\lambda(dy)
\end{equation}
where $c \geq 0$ and $\lambda$ is a Borel measure on $(0, \infty)$ satisfying $\int_{(0, \infty)}(1 \wedge y)\lambda(dy) < \infty$. We denote the law of $T(t)$ by $\rho_{t}$, for each $t \geq 0$. It is well-known that (see e.g. \cite{AG, App}) that if $\phi$ is a L\'{e}vy process in $G$ and $(T(t), t \geq 0)$ is an independent subordinator with Bernstein function $h$ then the process $\phi_{T} = (\phi_{T}(t), t \geq 0)$ is again a  L\'{e}vy process in $G$ where $\phi_{T}(t):= \phi(T(t))$ for all $t \geq 0$. Let $(P_{t}, t \geq 0)$ be the usual semigroup of convolution operators on $C(G)$ that is associated to the process $\phi$ and $(P^{h}_{t}, t \geq 0)$ be the corresponding semigroup of convolution operators associated to the process $\phi_{T}$. Then it is well-known (and easily verified) that
\begin{equation} \label{subsemi}
P^{h}_{t}f(\sigma) = \int_{(0, \infty)} P_{s}f(\sigma)\rho_{t}(ds)
\end{equation}
for all $t \geq 0, f \in C(G), \sigma \in G$.
In the sequel we will always take $\phi$ to be Brownian motion on $G$ as given by \eqref{sde}.
The following is implicit in \cite{App} but we give a proof for completeness.

\begin{prop} \label{subsemiFT} If $(P_{t}, t \geq 0)$ is the heat semigroup on $G$ then for all $f \in C(G), \pi \in \wG, t \geq 0$,
$$\widehat{P^{h}_{t}f}(\pi) = e^{-t h(\kappa_{\pi})}\widehat{f}(\pi).$$
\end{prop}

{\it Proof.} Using Fubini's theorem we have from (\ref{subsemi}) and Proposition \ref{ftheat}
\bean \widehat{P^{h}_{t}f}(\pi) & = & \int_{G}\int_{(0, \infty)} P_{s}f(\sigma)\rho_{t}(ds)\pi(\sigma^{-1})d\sigma \\
& = & \int_{(0,\infty)}\left(\int_{G}P_{s}f(\sigma)\pi(\sigma^{-1})d\sigma\right)\rho_{t}(ds)\\
& = & \int_{(0,\infty)}e^{-s\kappa_{\pi}}\rho_{t}(ds)\widehat{f}(\pi)\\
& = & e^{-t h(\kappa_{\pi})}\widehat{f}(\pi).~~~~~~~~\hfill \Box \eean

We now take $A=0$ in (\ref{identity}) and take the L\'evy process to be of the form $\phi_{T}$ as just described. For simplicity we also assume that $\psi$ only depends on the jumps of the process and so we write $\psi(\tau): = \psi(\cdot, \cdot, \tau)$ for each $\tau \in G$.  As before, we assume that $\psi$ is regular and that $\|\psi\|\leq 1$. In this case we have

\begin{eqnarray} \label{notime1}
\int_{G}S_{\psi}f(\sigma)g(\sigma)d\sigma &=&
 \int_{0}^{\infty}\int_{G}\int_{G}[(P^{h}_{s}f)(\sigma\tau)
- (P^{h}_{s}f)(\sigma)]\\
&\times& [(P^{h}_{s}g)(\sigma \tau)
- (P^{h}_{s}g)(\sigma)]\psi(\tau)\nu(d\tau)d\sigma ds\nonumber\end{eqnarray}
Now using (\ref{Planch1}) and Proposition \ref{subsemiFT} we obtain

\bean \int_{G}S_{\psi}f(\sigma)g(\sigma)d\sigma
 &=&  \int_{0}^{\infty}\int_{G}\sum_{\pi \in \wG}d_{\pi}e^{-2s h(\kappa_{\pi})}\tr[(\pi(\tau) - I_{\pi})\widehat{f}(\pi)\widehat{g}(\pi)(\pi(\tau)^{*} - I_{\pi})]\\
 &\times&\psi(\tau)\mu(d\tau)ds\\
& = & \int_{G}\sum_{\pi \in \wG}d_{\pi}\frac{1}{2 h(\kappa_{\pi})}\tr[(2I_{\pi} - \pi(\tau) - \pi(\tau)^{*})\widehat{f}(\pi)\widehat{g}(\pi)]\psi(\tau)\nu(d\tau) \eean
and so $S_{\psi}$ is a Fourier multiplier with
\begin{equation}\label{subordination-muilt}
m_{S_{\psi}}(\pi) = \frac{1}{2 h(\kappa_{\pi})}\int_{G}(2I_{\pi} - \pi(\tau) - \pi(\tau)^{*})\psi(\tau)\nu(d\tau),
\end{equation}
for $\pi \in \wG$.  From Theorem \ref{projections1} we obtain

\begin{corollary} Let $S_{\psi}$ be the operator with Fourier multiplier given by \eqref{subordination-muilt}.  \begin{enumerate}
\item If $\|\psi\|\leq 1$, then
\begin{equation}\label{subordination-1}
\|S_{\psi}f\|_p\leq (p^*-1)\|f\|_p,
\end{equation}
for all $1<p<\infty$, $f\in L^p(G)$.

\item If $\psi:G\to [b, B]$, with $-\infty<b<B <\infty$, then
\begin{equation}\label{subordination-2}
\|S_{\psi}f\|_p\leq C_{p, b, B}\|f\|_p,
\end{equation}
for all $1<p<\infty$, $f\in L^p(G)$, where $C_{p, b, B}$ is the constant in \eqref{martin}.
\end{enumerate}
\end{corollary}

To compare with results in \cite{BanBog} (see e.g. equation (32) therein), suppose that $G$ is connected as well as compact. Then $\exp: \g \rightarrow G$ is onto and so for every  $\sigma \in G, \sigma = \exp(X)$ for some $X \in \g$. Then $\pi(\sigma) = e^{i(-id\pi(X))}$ where we note that the matrix $-id\pi(X)$ is hermitian. Then we can write $2I_{\pi} - \pi(\sigma) - \pi(\sigma)^{*} = 2(I_{\pi}-\cos(id\pi(X))$. Let $\tilde{\nu}: = \nu \circ \exp$ so $\tilde {\nu}$ is a Borel measure on $\g$ and define $\tilde{\psi}:\g \rightarrow \R$ by $\tilde{\psi}:=\psi\circ\exp$. Then in this case we obtain
$$m_{S_{\psi}}(\pi) = \frac{1}{h(\kappa_{\pi})}\int_{\g}(I_{\pi}-\cosh(d\pi(X)))\tilde{\psi}(X)\tilde{\nu}(dX).$$

\subsection{Multipliers associated to central L\'{e}vy processes}

We say that a L\'{e}vy process $(\phi(t), t \geq 0)$ is {\it central} (or {\it conjugate invariant}) if its law $p_{t}$ is a central measure for all $t \geq 0$, i.e. for all $A \in {\calB}(G), \sigma \in G$
$$ p_{t}(\sigma A \sigma^{-1}) = p_{t}(A).$$ Such processes have been investigated in \cite{Liao1} and \cite{App}. It follows that the L\'{e}vy measure $\nu$ is also a central measure (see \cite{Liao1} p.1567.)
Recall that if $\mu$ is a Borel probability measure on a compact Lie group $G$ then its Fourier transform (or characteristic function\footnote{Note that in this paper we adopt the analysts convention when defining the Fourier transform of a function and the probabilists convention for that of a measure. We emphasize that no mathematical inconsistencies arise from these choices.}) is defined by $\widehat{\mu}(\pi):= \int_{G}\pi(\sigma)\mu(d\sigma)$ for all $\pi \in \wG$. Moreover the matrices $\{\widehat{\mu}(\pi), \pi \in \wG\}$ uniquely determine $\mu$ (see e.g. Theorem 1 in \cite{KI}.) Returning to L\'{e}vy processes, it is easy to verify that for all $t \geq 0, \pi \in \wG, 1 \leq i,j \leq d_{\pi}$ we have
\begin{equation} \label{semFT}
\widehat{p_{t}}(\pi)_{ij} = P_{t}\pi_{ij}(e).
\end{equation}

Now suppose that the L\'{e}vy process $(\phi(t), t \geq 0)$ is central. By a slight generalization of the argument of Proposition 2.1 in \cite{App} we deduce that for all $\pi \in \wG$, there exists $\alpha_{\pi} \in \Co$ with $\Re(\alpha_{\pi}) \leq 0$ so that
\begin{equation} \label{centFT}
\widehat{p_{t}}(\pi) = e^{t \alpha_{\pi}}I_{\pi}.
\end{equation}
It follows easily from this that for all $t \geq 0, \pi \in \wG, f \in C(G)$
\begin{equation} \label{centFTsemi}
\widehat{P_{t}f}(\pi) = e^{t \alpha_{\pi}}\widehat{f}(\pi)
\end{equation}
We now apply \eqref{identity}. For simplicity we again consider the case where $A$ is a constant matrix and $\psi=\psi(\tau)$ depends only on the jumps of the L\'{e}vy process. In this case we can apply Plancherel's theorem \eqref{Planch1} using very similar arguments to those applied in sections 4.1 and 4.4 to deduce that $S_{A, \psi}$ is  Fourier multiplier. Indeed for all $f \in C^{\infty}(G), \pi \in \wG$ we have
$$ \widehat{S_{A, \psi}f}(\pi) = m_{A,\psi}(\pi){\widehat f}(\pi),$$ where
\bean m_{A,\psi}(\pi) & = &  \frac{1}{\Re(\alpha_{\pi})}\sum_{i,j=1}^{n}A_{ji}d\pi(X_{i})d\pi(X_{j})\\
& - & \frac{1}{2\Re(\alpha_{\pi})}\int_{G}(2I_{\pi} - \pi(\tau) - \pi(\tau)^{*})\psi(\tau)\nu(d\tau).
\eean
Two of our previously considered examples are in fact special cases of this formula. Indeed for the second order Riesz transform we have $\psi = 0$ and $\alpha_{\pi} = -\kappa_{\pi}$ while for subordinated Brownian motion, $A = 0$ and $\alpha_{\pi} = -h(\kappa_{\pi})$.

Under further constraints, we can find an explicit general formula for $\alpha_{\pi}$. The following is essentially contained in Theorem 4 (c) p.1567 in \cite{Liao1} under some stronger assumptions. It gives a L\'{e}vy-Khintchine type formula for the Fourier transform of the law of a central L\'{e}vy process. We will need the group character  $\chi_{\pi}(\cdot):= \tr(\pi(\cdot))$ and the normalized character $\varrho_{\pi}(\cdot):=\frac{1}{d_{\pi}}\chi_{\pi}$ for each $\pi \in \wG$.

\begin{theorem}
Assume that $(\phi(t), t \geq 0)$ is a central L\'{e}vy process on a compact, connected semi-simple Lie group $G$. Then for all $\pi \in \wG$,
$$ \alpha_{\pi} = -c \kappa_{\pi} + \int_{G}(\varrho_{\pi}(\tau) - 1)\nu(d\tau),$$
where $c \geq 0$.

\end{theorem}

{\it Proof.} First assume that the L\'{e}vy measure $\nu$ has a finite first moment, i.e. $\int_{G}|x_{i}(\sigma)|\nu(d\sigma) < \infty$ for all $1 \leq i \leq n$. It follows from Propositions 4.4 and 4.5 in \cite{Liao} p.99 that for all $f \in C^{2}(G), \sigma \in G$,
$$ {\calL}f(\sigma) = c \Delta f(\sigma) + \int_{G}(f(\sigma\tau)- f(\sigma))\nu(d\tau).$$
For each $\pi \in \wG$, define the matrix ${\calL}(\pi)$ by ${\calL}(\pi)_{ij}:=({\calL}\pi_{ij})(e)$ for $1 \leq i,j \leq d_{\pi}$. Then
$${\calL}(\pi) = - c \kappa_{\pi}I_{\pi} + \int_{G}(\pi(\tau) - I_{\pi})\nu(d\tau),$$ and hence by \eqref{semFT},
$ \widehat{p_{t}}(\pi) = e^{t {\calL}(\pi)}$ for all $t \geq 0$. Comparing this identity with \eqref{centFT}, we deduce that
$$ \alpha_{\pi}I_{\pi} = -c \kappa_{\pi}I_{\pi} + \int_{G}(\pi(\tau) - I_{\pi})\nu(d\tau),$$
and we obtain the required result on taking the trace (in $V_{\pi}$) of both sides of this last equation.
Now we turn to the general case and take $(\phi(t), t \geq 0)$ to be an arbitrary central L\'{e}vy process with characteristics $(b, a, \nu)$ and infinitesimal generator ${\calL}$. Again by Proposition 4.4 in \cite{Liao}, we must have $a = cI$ for some $c \geq 0$ and $b=0$ by the argument in \cite{Liao}, Proposition 4.5. So the process has characteristics $(0, cI, \nu)$.  Let $(U_{m}, m \in \mathbb{N})$ be a sequence of Borel sets in $G$ so that $U_{m} \uparrow G$ as $m \rightarrow \infty$  and  $\int_{U_{m}}|x_{i}(\sigma)|\nu(d\sigma) < \infty$ for all $1 \leq i \leq n$. To see that such sets can always be constructed, recall that $G$ is equipped with a bi-invariant metric for which the mappings $\tau \rightarrow \sigma \tau \sigma^{-1}$ (where $\tau \in G$) are isometries for all $\sigma \in G$. So we may take e.g. $U_{m} = B_{\frac{1}{m}}(e)^{c}$. Now define the central L\'{e}vy measure  $\nu_{m}$ on $G$ by $\nu_{m}(A): = \nu(A \cap U_{m})$ for all $A \in \calB(G),m \in \mathbb{N}$. For each $m \in \mathbb{N}$, let $(\phi_{m}(t), t \geq 0)$ be a central L\'{e}vy process with characteristics $(0, cI, \nu_{m})$ and infinitesimal generator ${\calL}_{m}$. Then it is easily verified that for all $f \in C^{2}(G)$,
$$ \lim_{m \rightarrow \infty}{\calL}_{m}f = {\calL}f.$$
If $p_{t}^{(m)}$ is the law of $\phi_{m}(t)$ then for all $\pi \in \wG$,
$$\widehat{p_{t}^{(m)}}(\pi) = \exp{\left\{ -c \kappa_{\pi} + \int_{U_{m}}(\varrho_{\pi}(\tau) - 1)\nu(d\tau)\right\}},$$ but from the construction in the first part of the proof, we have $\lim_{m \rightarrow \infty}\widehat{p_{t}^{(m)}}(\pi) = \widehat{p_{t}}(\pi)$ and the result follows. $\hfill \Box$

\section{$\bR^n$, revisited}

In this section we apply the results of \S3 to Euclidean  space $\bR^n$ and compute the Fourier multipliers for the operators $S_{A,\psi}$ explicitly under some assumptions on $A$ and  $\psi$.  When $A$ is constant and $\psi$ is just a function of the jumps this was already done in \cite{BanBieBog} and \cite{BanBog}.  Nevertheless, the formulas derived in \S3 elucidate these results further and provide a more uniform approach to these L\'evy-Fourier multipliers.

We start by recalling various standard notations and facts about Fourier transforms and L\'evy processes in $\bR^n$. We use the following normalization for the Fourier transform $\widehat{f}$ of $f \in {\calS}(\R^{n})$ where ${\calS}(\R^{n})$ is the usual Schwartz space of rapidly decreasing functions on $\R^{n}$. For all $\xi \in \R^{n}$,
\begin{equation}\label{Fo}
\widehat{f}(\xi)=\int_{\bR^n} e^{2\pi i\xi\cdot x}f(x)dx
\end{equation}
and for all $x \in \R^{n}$,
\begin{equation}\label{Foin}
f(x)=\int_{\bR^n} e^{-2\pi ix\cdot\xi}\widehat{f}(\xi)d\xi
\end{equation}
so that Plancherel's identity takes the form

\begin{equation}\label{plancherel}
\int_{\bR^n} f(x)g(x)dx=\int_{\bR^n}\widehat{f}(\xi)\overline{\widehat{g}(\xi)} d\xi,
\end{equation}
for $f, g \in {\calS}(\R^{n})$.
With this normalization, if $\nabla=\left(\frac{\partial}{\partial x_1}, \dots, \frac{\partial}{\partial x_n}\right)=\left(\partial_1, \dots, \partial_n\right)$ denotes the standard gradient in $\bR^n$,  then
\begin{equation} \label{nabla}
\widehat{\partial_j f}(\xi)=(-2\pi i) \xi_j \hat{f}(\xi) \,\,\,\,\text{and}\,\,\,\, \widehat{\nabla f}(\xi)=(-2\pi i)\xi \hat{f}(\xi).
\end{equation}

Recall that a Borel measure  $\nu$
on $\bR^n$ with $ \nu(\{0\})=0$ and
 \begin{equation}\label{levymu1}
 \int_{\bR^n} \frac{|y|^2}{1+|y|^2} \,\nu(dy)<\infty
\end{equation}
is called  a  L\'evy measure.  (Note that the definition of the L\'evy measure in $\bR^n$ coincides with that given by \eqref{levy-Lie} in \S2.2 for general Lie groups.)
We denote a L\'evy process in $\bR^n$ by $(X(t), t\geq 0)$.  The celebrated  L\'evy-Khintchine formula
\cite{App1}
guarantees the existence of a triple
$\left( b, a, \nu \right)$
such that the characteristic function of the process is given by
$\E\left[\, e^{i \xi \cdot X(t)}\,\right]  = e^{t
\rho(\xi)},$
where
\begin{equation}\label{levykin}
\rho(\xi)= i b\cdot \xi  - a\xi\cdot \xi +
\int_{\bR^n}\left[\,  e^{i\,\xi \cdot y}- 1 - i (\xi\cdot y) \,\bI_{B_{1}(0)} (y) \,\right]
\nu(dy).
\end{equation}
Here,  $b=(b_1, \dots, b_n) \in \bR^n$,   $a=(a_{ij})$ is a non-negative $n \times n$ symmetric matrix,
$\bI_{B_1(0)}$
is the indicator function of the unit ball $B_{1}(0)\subset\bR^n$ centered at the origin and $\nu$
is a L\'evy measure  on $\bR^n$. We remark that the L\'evy-Khintchine formula may be deduced as a corollary of Hunt's theorem (Theorem \ref{Hunt}) in the case where the group $G$ is $\R^{n}$ (see \cite{Hunt}, \S6, pp.281-3.)
 We decompose $\rho$ into its real and imaginary parts so that
\begin{equation}\label{real-Imag}
\rho(\xi)=\Re{\rho}(\xi)+i\Im{\rho}(\xi),
\end{equation}
where
\begin{equation}\label{real}
\Re{\rho}(\xi)=- a\xi\cdot \xi +
\int_{\bR^n}\left[\,  \cos{(\xi \cdot y)}- 1\,\right]
\nu(dy)
\end{equation}
and
\begin{equation}\label{imag}
\Im{\rho}(\xi)=b\cdot \xi+ \int_{\bR^n}\left[\, \sin{(\xi \cdot y)}-  (\xi\cdot y) \,\bI_{B_{1}(0)} (y) \,\right]\,
\nu(dy).
\end{equation}
We note that the convergence of the integrals in \eqref{real} and \eqref{imag} follows immediately from \eqref{levymu1}.  Also, observe that $\Re{\rho}(-\xi)=\Re{\rho}(\xi)$ and that $\Re{\rho}(\xi)\leq 0$.

With this notation, the semigroup of the L\'evy process $(X(t), t\geq 0)$  acting on $f\in {\calS}(\bR^n)$ is given as a pseudo-differential operator by
\begin{eqnarray}\label{semigroup}
P_tf(x)=\E[f(X(t)+x)] &=& \E\left(\int_{\bR^n}e^{-2\pi i(X(t)+x)\cdot\xi}\,\widehat{f}(\xi)d\xi\right)\nonumber\\
&=&\int_{\bR^n}e^{t\rho(-2\pi\xi)}\,e^{-2\pi ix\cdot\xi}\widehat{f}(\xi)d\xi.
\end{eqnarray}
Since ${\calS}(\R^{n}) \subseteq \Dom({\calL})$ we also obtain the pseudo-differential operator representation
\begin{equation}
{\calL}f(x)=\int_{\bR^n} \rho(-2\pi\xi) e^{-2\pi ix\cdot\xi}\,\widehat{f}(\xi)d\xi.
\end{equation}
Thus
\begin{equation}\label{fouriersemigroup}
\widehat{P_tf}(\xi)=e^{t\rho(-2\pi\xi)}\,\widehat{f}(\xi).
\end{equation}
Also, using \eqref{nabla} we can write $\mathcal{L}$ as
 \begin{eqnarray*}
\mathcal{L}f(x)=  b^i\partial_i f(x)&+& a^{ij}\partial_i\partial_j f(x)\\
&+& \int\Big[f(x+y)-f(x)-y\cdot\nabla f(x)\bI_{B_{1}(0)} (y)\Big] \nu(dy).
\end{eqnarray*}
A detailed account of these results may be found in \cite{App1} \S3.2.2, pp.163-9.

With $\Lambda$ such that $\Lambda \Lambda^{T} = 2a$, we set, as in \S3, $\nabla_{Y}=\Lambda\nabla$.  With this notation we have the following corollary of Theorem \ref{projections1}.

\begin{corollary}\label{projection-Rn}
Suppose $A \in L^{\infty}(\R^{+} \times \bR^n, \mbox{Leb} \times m_{R};M_{n}(\R))$ and $\psi \in L^{\infty}(\R^{+} \times \bR^n \times \bR^n; \bR)$ with $(A, \psi)$ regular and $||A|| \vee ||\psi|| \leq 1$.  There exists a bounded linear operator $S_{A, \psi}$ on $L^{p}(\bR^n)$ where  $1<p<\infty$ for which

\begin{eqnarray} \label{identity-Rn}
\int_{\bR^n}S_{A,\psi}f(x)g(x)dx& = & \int_{0}^{\infty}\int_{\bR^n} A(s, x)\nabla_{Y}(P_{s}f)(x)\cdot \nabla_{Y}(P_{s}g)(x)dx ds \nonumber \\
& + & \int_{0}^{\infty}\int_{\bR^n}\int_{\bR^n}[(P_{s}f)(x+y)
- (P_{s}f)(x)][(P_{s}g)(x+y)
- (P_{s}g)(x)]\\
&\times& \psi(s,x, y)\nu(dy)dx ds,\nonumber
\end{eqnarray}
for al $f, g \in C_{c}^{\infty}(G)$.
Furthermore, for all $f\in L^p(\bR^n)$ and $g\in L^q(\bR^n)$, $\frac{1}{p}+\frac{1}{q}=1$,
\begin{equation}\label{bound1-Rn}
\Big|\int_{\bR^n}S_{A,\psi}f(x)g(x)dx \Big| \leq (p^*-1)\|f\|_p\,\|g\|_q
\end{equation}
and
\begin{equation}\label{bound2-Rn}
\|S_{A,\psi}f\|_p\leq (p^*-1)\|f\|_p.
\end{equation}
\end{corollary}

Let us now suppose that $A=A(s)$ is only a function of $s$ (time) and that $\psi(s, y)$ is only a function of $s$ (time) and $y$ (jumps).  Set

\begin{equation}
I=\int_{0}^{\infty}\int_{\bR^n} A(s)\nabla_{Y}(P_{s}f)(x)\cdot \nabla_{Y}(P_{s}g)(x)dx ds
\end{equation}
and
\begin{equation}
II= \int_{0}^{\infty}\int_{\bR^n}\int_{\bR^n}[(P_{s}f)(x+y)
 (P_{s}f)(x)][(P_{s}g)(x+y)
- (P_{s}g)(x)]\psi(s, y)\nu(dy)dx ds.
\end{equation}
Applying Plancherel's identity, \eqref{plancherel}, we find that
\begin{eqnarray}\label{I}
I &=&\int_{\bR^n}\left\{4\pi^2\int_0^{\infty} [A(s)\Lambda\xi\cdot \Lambda\xi]e^{2s\Re{\rho(2\pi\xi)}}ds\right\} \hat{f}(\xi)\overline{\hat{g}}(\xi)d\xi\\
&=&\int_{\bR^n} m_1(\xi)\hat{f}(\xi)\overline{\hat{g}}(\xi)d\xi,\nonumber
\end{eqnarray}
where
\begin{equation}\label{mult-1}
m_1(\xi)=4\pi^2\int_0^{\infty} [A(s)\Lambda\xi\cdot \Lambda\xi]e^{2s\Re{\rho(2\pi\xi)}}ds.
\end{equation}

In the same way,
\begin{eqnarray}\label{II}
II&=&\int_{\bR^n}\left\{\int_{\bR^n}\int_0^{\infty}e^{2s\Re{\rho}(2\pi\xi)}|e^{-2\pi\,i\xi\cdot y}-1|^2 \psi(s, y)ds d\nu(y)  \right\} \hat{f}(\xi)\overline{\hat{g}}(\xi)d\xi\nonumber\\
&=&\int_{\bR^n}\left\{2\int_{\bR^n}\int_0^{\infty}e^{2s\Re{\rho}(2\pi\xi)}\left(1-\cos(2\pi\xi\cdot y)\right)\psi(s, y)ds \nu(dy)   \right\} \hat{f}(\xi)\overline{\hat{g}}(\xi)d\xi\\
&=&\int_{\bR^n} m_2(\xi)\hat{f}(\xi)\overline{\hat{g}}(\xi)d\xi,\nonumber
\end{eqnarray}
where
\begin{equation}\label{multiplier-2}
m_2(\xi)=2\int_{\bR^n}\int_0^{\infty}e^{2s\Re{\rho}(2\pi\xi)}\left(1-\cos(2\pi\xi\cdot y)\right)\psi(s, y)ds d\nu(y).
\end{equation}

Thus under the assumption that  $A=A(s)$ and $\psi = \psi(s,y)$ we conclude that the operator $S_{A,\psi}$  is a Fourier multiplier with
\begin{equation}
\widehat{S_{A,\psi}f}(\xi) =m(\xi) \widehat{f}(\xi),
\end{equation}
where
\begin{eqnarray}\label{multiplier-3}
m(\xi)&=&4\pi^2\int_0^{\infty} [A(s)\Lambda\xi\cdot \Lambda\xi]e^{2s\Re{\rho(2\pi\xi)}}ds\nonumber\\
&+&2\int_{\bR^n}\int_0^{\infty}e^{2s\Re{\rho}(2\pi\xi)}\left(1-\cos(2\pi\xi\cdot y)\right)\psi(s, y)ds \nu(dy)
\end{eqnarray}

Recalling from  \eqref{real} that
$$
-2\Re{\rho}(2\pi \xi)= 8\pi^2a\xi\cdot \xi +
2\int_{\bR^n}\left(1- \cos{(2\pi\xi \cdot y)}\right) \nu(dy)
$$
and using the fact that $2a\xi\cdot \xi=\Lambda\xi\cdot \Lambda\xi$,   we see that under the assumption  $||A|| \vee ||\psi|| \leq 1$, $m\in L^{\infty}(\bR^n)$ and $\|m\|_{\infty}\leq 1$.  Furthermore, if the matrix $A$ has constant entries and the function $\psi$ only depends on $y$,  then a simple computation gives that

\begin{eqnarray}\label{multiplier-4}
m(\xi)=\frac{4\pi^2 A\Lambda\xi\cdot \Lambda\xi +
2\int_{\bR^n}\left(1- \cos{(2\pi\xi \cdot y)}\right)\psi(y) \nu(dy)}{8\pi^2a\xi\cdot \xi +
2\int_{\bR^n}\left(1- \cos{(2\pi\xi \cdot y)}\right) \nu(dy)}.
\end{eqnarray}
Note that $m(\xi)=\widetilde{m}(2\pi\xi)$,
where
$$
\widetilde{m}(\xi)=\frac{\frac{1}{2}\left(\Lambda^TA\Lambda\right)\xi\cdot \xi +
\int_{\bR^n}\left(1- \cos{(\xi \cdot y)}\right)\psi(y) \nu(dy)}{a\xi\cdot \xi +
\int_{\bR^n}\left(1- \cos{(\xi \cdot y)}\right) \nu(dy)}.
$$

From the Remark \ref{hilbert-valued} we know that the above multipliers have the same bounds if the matrix has complex entires and the function $\psi$ takes values in the complex plane.  We summarize these results in the following
\begin{theorem}\label{multipliersRn}
Let $A\in M_{n}(\Co)$ and $\psi:\bR^n\to\Co$ be such that $\|A\|\vee \|\psi\|\leq 1$. Suppose $a\in M_{n}(\bR)$ is symmetric and non-negative definite.  Let $\Lambda\in M_{n}(\bR)$ be such that $\Lambda\Lambda^T=2a$.  Then the $L^{\infty}$ function
\begin{equation}\label{multiplier-5}
{m}(\xi)=\frac{\frac{1}{2}\left(\Lambda^TA\Lambda\right)\xi\cdot \xi +
\int_{\bR^n}\left(1- \cos{(\xi \cdot y)}\right)\psi(y) \nu(dy)}{a\xi\cdot \xi +
\int_{\bR^n}\left(1- \cos{(\xi \cdot y)}\right) \nu(dy)},
\end{equation}
where $\nu$ is a L\'evy measure, defines an $L^p$-multiplier operator $S_{A,\psi}$  with \
\begin{equation}\label{bound2Rn}
\|S_{A,\psi}f\|_p \leq (p^*-1)\|f\|_p,
\end{equation}
for all $f\in L^p(\bR^n)$, $1<p<\infty$.    The inequality is sharp.
\end{theorem}
This formula gives  the multipliers studied in \cite{BanBieBog} and \cite{BanBog}.  For several concrete examples of classical multipliers that arise from this formula, see \cite{Ban1} and \cite{BanBieBog}.  In particular, if we take $\nu=0$, $a=I$ and we write $A=(A_{jk})$, then
$$
m(\xi)=\frac{\sum_{j, k=1}^n A_{jk}\xi_j\xi_k}{|\xi|^2}.
$$
With the right choice of $A$, this gives the multiplier $m(\xi)=\frac{\xi_j^2-\xi_k^2}{|\xi|^2}$.   This corresponds to the second order Riesz transforms $R_k^2-R_j^2$ which by \cite{GMS} (see also \cite{BanOse}) has norm $(p^*-1)$.  Hence the bound is sharp.

\end{document}